\documentclass[12pt,a4paper]{article}
\usepackage{amsmath,amsfonts} 

\usepackage{graphicx}
\usepackage{amsthm}

\newtheorem{theorem}{Theorem}

\newtheorem{definition}[theorem]{Definition}

\newtheorem{proposition}[theorem]{Proposition}
\newtheorem{remark}[theorem]{Remark}

\newcommand{\SO}{\mathcal{O}}
\newcommand{\fg}{\mathfrak{g}}
\newcommand{\fh}{\mathfrak{h}}
\newcommand{\fl}{\mathfrak{l}}
\newcommand{\fz}{\mathfrak{z}}
\newcommand{\ft}{\mathfrak{t}}
\newcommand{\rk}{\operatorname{rank}}
\newcommand{\CC}{\mathbb{C}}
\newcommand{\SP}{\mathcal{P}}

\begin{document}

\begin{center}
{\LARGE Projective moduli space of semistable principal}

{\LARGE sheaves for a reductive group
\footnote{
Article presented by I.S. in the conference held at Catania 
(11-13 April 2001), dedicated to the 
60th birthday of Silvio Greco.
The proceedings of this conference will appear in
Le Matematiche (Catania).\\
Both authors are members of EAGER (EC FP5 Contract no.
HPRN-CT-2000-00099). T.G. was supported by a postdoctoral
fellowship of Ministerio de Educaci\'on y Cultura (Spain).
} 
}

\bigskip

\bigskip

{\large Tom\'{a}s L. G\'{o}mez and Ignacio Sols}\bigskip \bigskip

\bigskip

\bigskip

\bigskip

\end{center}

\hspace{5cm}\textit{In the sixtieth anniversary
of Silvio Greco, }

\hspace{5cm}\textit{to Nadia and Silvio.}

\bigskip

\section{\protect\bigskip Introduction}

\noindent This contribution to the homage to Silvio Greco is mainly an
announcement of results to appear somewhere in full extent, explaining their
development from our previous article \cite{G-S1} on conic bundles.

In \cite{N-S} and \cite{S1} Narasimhan and Seshadri defined stable 
bundles on a curve
and provided by the techniques of Geometric Invariant Theory (GIT) developed
by Mumford \cite{Mu} a projective moduli space of the stable equivalence 
classes of
semistable bundles. Then Gieseker \cite{Gi} and Maruyama 
\cite{MaI} \cite{MaII} 
generalized this
construction to the case of a higher-dimensional projective variety,
obtaining again a projective moduli space 
by also allowing torsion-free sheaves.
Ramanathan \cite{R1} \cite{R2} has 
provided the moduli space of semistable principal
bundles on a connected reductive group $G$, thus generalizing the Narasimhan
and Seshadri notion and construction, which then becomes the
particular case 
$G=Gl (n,\mathbb{C).}$ 

Faltings \cite{Fa} has considered the moduli stack of 
principal bundles
on semistable curves. For $G$ orthogonal or symplectic he
considers a torsion-free sheaf with a quadratic form, and
he also defines a notion of stability. For general
reductive group $G$ he uses the approach of loop groups.
Sorger \cite{So} had considered a similar problem. 
He works on a curve $C$ (not necessarily smooth) on a smooth
surface $S$, and  constructs
the moduli space of torsion free sheaves on $C$ together with
a symmetric form taking values on the dualizing sheaf $\omega_C$.

In the talk ``open problems on principal bundles''
closing the conference on ``vector bundles on algebraic curves 
and Brill-Noether theory'' at Bad Honeff
2000, prof. Narasimhan proposed the problem of generalizing the work 
of the late
Ramanathan to the case of higher-dimensional varieties and to the case of
positive characteristics. We solve the first problem by providing a suitable
definition of principal sheaf on a higher-dimensional projective variety $X$
over the complex field, and a definition of their (semi)stability, which in
case $\dim X=1$ is that of Ramanathan, and for which a projective moduli
space can be obtained.

We start by recalling our notion of (semi)stable conic bundles, i.e.
symmetric $(2,0)$-tensors $\varphi :E\otimes E\rightarrow \SO _{C}$ of 
$\rk E=3$ on an algebraic curve $C$, and their projective moduli
space, 
notion and moduli space which have been generalized to the case of $(s,0)$ 
tensors on a
curve by Schmitt \cite{Sch} with the purpose of dealing with (semi)stable 
objects 
$\varphi :E^{\rho}\rightarrow M$,
where $E^\rho$ is the vector bundle associated to a vector 
bundle $E$ and an
arbitrary representation $\rho$ of $G=Gl (n,\mathbb{C})$,
and $M$ is a line bundle. In case the
symmetric $(2,0)$ tensor is of maximal rank at all points
and $\det(E)\cong \SO_C$, 
i.e. the case when 
$(E,\varphi )$ is just a principal $SO(3,\mathbb{C})$-bundle, our notion of
(semi)stability is drastically simplified and becomes 
equivalent to Ramanathan's
notion of (semi)stability. We then generalize to higher dimension,
with techniques of Simpson \cite{Si} and
Huybrechts-Lehn \cite{H-L},
the notion
and coarse projective moduli space of (semi)stable $(s,0)$-tensors, 
by allowing $E$
to be a torsion free sheaf and those symmetric or antisymmetric and nowhere
degenerate provide thus the moduli space 
of principal sheaves on $G=O(n,\mathbb{C)}$, 
$Sp(n,\mathbb{C),}$ $SO(2n+1,\mathbb{C})$, 
the remaining classical group $SO(2n,\mathbb{C})$ 
requiring a special treatment which fortunately does not alter the
notion of (semi)stability.

Then we cope with the problem of an arbitrary connected reductive group $G$,
by defining principal sheaves as $(2,1)$ tensors, i.e. torsion free 
sheaves $E$ and $\varphi :E\otimes E\rightarrow E^{\ast\ast}$, which on 
the points of the open
set $U_{E}$ where $E$ is locally free are isomorphic to the structure 
tensor 
$\varphi _{\fg'}:\fg'\otimes \fg'\rightarrow \fg' $ 
of the Lie algebra $\fg'$ tangent to the commutator 
$G^{\prime }=[G,G],$
together with a $G\rightarrow Aut(\fg')$ reduction of the associated
principal bundle on $U_{E}$. The (semi)stability is defined as just the
one of the $(2,1)$ tensor $(E,\varphi )$ and leads to a coarse projective
moduli space, reducing to Ramanathan's (semi)stability and moduli
space in case $\dim X=1 $.

This announcement note consists mainly of the precise definitions and
statements of such objects and results.

\section{Conic bundles.}

Let $X$ be a complete, smooth, connected curve, and fix a positive 
rational $\tau >0$. A conic bundle on $X$ of degree $d$ is a rank $3$ 
symmetric $(2,0)$-tensor on $X$, i.e. a vector bundle $E$ 
of rank $3$ and degree $d$, together
with a nonzero homomorphism 
\[
\varphi :E^{2}=Sym^{2}E\rightarrow L
\]
Where $L$ is a line bundle.
We say it is (semi)stable if

1) For any subbundle $F\subseteq E$, it is 
\[
\frac{\deg F-c_{\varphi }(F)\tau }{\rk F}\;(\leq )\;
\frac{\deg  E-2\tau }{\rk E}
\]
where 
\[
c_{\varphi }(F)=\left\{ 
\begin{array}{l}
2,\;\text{if}\;\varphi (F^{2})\neq 0\quad
(\text{i.e.} F \;\text{not isotropic)}\\ 
1,\;\text{if}\;\varphi (F^{2})=0\;\operatorname{ and }\;\varphi (FE)\neq 0 \\ 
0, \;\text{if}\;\varphi (FE)=0 \quad\text{(i.e. }F\;{\rm singular)}
\end{array}
\right. 
\]

2) For all critical flags $F_{1}\subseteq F_{2}\subseteq E$, i.e. $F_{1}$
of rank $1$, $F_{2}$ of rank $2$, $\varphi(F_1F)\neq0$,
$\varphi(F_2F_2)\neq 0$ and $\varphi
(F_{1}F_{2})=0$ (at a general point of $X$ these are a point of the conic
and its tangent line) , the following inequality holds 
\[
\deg F_{1}+\deg F_{2}\;(\leq )\; \deg E
\]
(By the expression (semi)stable we always mean both semistable and stable,
and then by the symbol $(\leq )$ we mean $\leq $ and $<$, respectively).
As usual, there is a notion of stable equivalence classes of
semistable objects (see \cite{G-S1} for the definition), and then 
it is proved in \cite{G-S1}, by the use of GIT, the following

\begin{theorem}
\label{1}
There is a projective coarse moduli space of stable equivalence classes of
semistable conic bundles of degree $d$ and parameter $\tau$, 
on a smooth, complete, connected curve.
\end{theorem}

If $\det(E)\cong \SO_X$, $L\cong \SO_X$ and $\varphi$ is
nowhere degenerate,i.e. such that 
$\rk \varphi (x)=3$ for all $x\in X$, 
which amounts to a 
\textit{principal $SO(3)$-bundle on $X$}, then
the condition \ 2), independent of the
parameter $\tau $ is enough for the definition of
(semi)stability, thus leading to a projective coarse moduli space as in 
Theorem \ref{1}.

Recall that in \cite{R1}, \cite{R2}, a definition of (semi)stable principal bundle $P$
on a curve, for a connected, reductive group $G$ was already given: if for
all reduction $P(H)$ of $P$ to a maximal parabolic subgroup $H\subseteq G$,
the vector bundle $P(H,\fh)$ associated to $P$ by the adjoint representation
of $H$ in its tangent Lie algebra $\fh$, has 
\[
\deg P(H,\fh)\;(\leq )\; 0 
\]
In fact Ramanathan obtains in \cite{R2} a projective coarse moduli
space of stable
equivalence classes of semistable principal $G$-bundles of fixed topological
type and our result for $SO(3)$-bundles on $X$ becomes a particular case of
Ramanathan's result, because it is proved in \cite{G-S1} that condition 2 is
equivalent to the notion of Ramanathan.

Rank $2$ bundles correspond, after projectivization, to geometrically ruled
surfaces, and properties of the (semi)stable objects have been largely
studied since their definition in \cite{N-S} and \cite{S1}. Our definition of
(semi)stable conic bundles opens analogous problems. For instance we would
like to express here the following \textit{conjecture}. If has been proved
in \cite{C-S}, for a semistable scroll of $\mathbb{P}^{r}$ of degree $d$ and
irregularity $q$ which is special (i.e. $r$ distinct of 
Riemann-Roch number $d+1-2q$), the existence of a hyperplane 
containing $r-1$ lines of the
ruling, which amounts to the upper bound $d-(r-1)$ for the degree of a
unisecant curve of the ruled surface, a problem posed by Severi in 
\cite{Se} (the
analogous bound being trivial in the nonsemistable case). Most probably, for
a special semistable conic bundle of $\mathbb{P}^{r}$ there is a hyperplane
containing $\left[ \frac{r-2}{2}\right] $ of its conics, thus leading to an
analogous upper bound of the minimal degree of a bisecant curve of the
surface (and so on).

\section{Principal sheaves for a classical group}
\label{sec3}

Let $X$ be a smooth, projective complex variety of dimension $n$.

\begin{definition}
A tensor field, or just a tensor, on $X,$ is a pair $(E,\varphi )$ 
consisting of a torsion free sheaf $E$ and an homomorphism 
\[
\varphi :\otimes^{s}E\rightarrow \SO _{X},
\] 
the rank and Chern classes of the tensor
being called those of $E$. Let $\sigma $ be a positive rational polynomial
of degree at most $n-1$ (i.e. rational coefficients, and positive leading
coefficient). The tensor is said to be $\delta$-(semi)stable if for all
weighted filtration (E.,m.) of $E$, i.e. subsheaves $E_{1}\subset ...\subset
E_{t}\subset E_{t+1}=E$ and positive integers $m_{1},...,m_{t}$, it is 
\[
\sum m_{i}(r\chi _{E_{i}}-r_{i}\chi _{E})+
\delta\;\mu (E.,m.,\varphi )(\leq )0
\]
where $r,r_{i},\chi _{E},\chi _{E_{i}}$ are the ranks and Hilbert
polynomials of $E,E_{i},$ and $\mu $ is defined as
\[
\mu =\min \left\{ \lambda _{i_{1}}+...+\lambda _{i_{s}}|\varphi
(E_{i_{1}}\otimes ...\otimes E_{i_{s}})\neq 0\right\} 
\]
where $\lambda _{1}<...<\lambda _{s}$  are integers with 
$\lambda _{i}-\lambda _{i-1}=m_{i}$ 
and 
\[
\sum \lambda _{i} \rk(E_{i}/E_{i-1})=0.
\]

\end{definition}

In \cite{G-S2} the definition is slightly more general, and we
prove the following

\begin{theorem}
There is a coarse projective moduli space 
of $\delta $-stable equivalence classes of 
$\delta$-semistable tensors on a projective variety $X$, of fixed Chern
classes and rank.
\end{theorem}

The proof has two parts: first, show that the family consisting of such
objects is bounded (remark that for $\delta$-semistable $(E,\varphi $), the
torsion free sheaf $E$ needs not be semistable). Second, proceed with
the techniques of Simpson \cite{Si} and Huybrechts-Lehn \cite{H-L}, 
starting by considering an integer $m\gg 0$ such that
all torsion free sheaves in the family are generated by global sections and
have $H^{0}(E(m))=\chi _{E}(m)$. For each member of the family choose an
isomorphism $\beta $ of $H^{0}(E(m))$ with a fixed complex vector space $V$
of dimension $\chi (E(m))$, thus obtaining a quotient 
\[
V\otimes \SO _{X}(-m)\simeq H^{0}(E(m))\otimes \SO
_{X}(-m)\longrightarrow E
\]
inducing, for $l$ high enough, a quotient 
\[
q:V\otimes H^{0}(\SO _{X}(l -m))\longrightarrow H^{0}(E(l )).
\]
Consider also the induced homomorphism 
\[
\psi :V^{\otimes s}\longrightarrow 
H^{0}(E(m)^{\otimes s})\longrightarrow
H^{0}(\SO_X(sm))
\]
We then obtain an element of 
\[
\mathbb{P}\Big(\bigwedge{}^{\chi^{}_E (l)}_{}\big( V^{\ast }\otimes 
H^{0}(\SO _{X}(l-m))^{\ast} \big)\Big)\times 
\mathbb{P}\Big(V^{\ast\otimes s }\otimes H^0(\SO_X(sm)\Big)
\]
which we consider included in projective space by the linear system 
$|\SO(n_{1},n_{2})|$ with 
\[
\frac{n_{2}}{n_{1}}=\frac{\chi_E (l )\delta (m)-
\delta (l )\chi_E (m)}{\chi_E (m)-s\delta (m)}
\]
This assignation embeds in a projective space $\mathbb{P}$
the scheme $R$ of
triples $(E,\varphi ,\beta )$, with $(E,\varphi )$ being a 
$\delta$-semistable tensor of the given rank and Chern classes and 
$\beta $ a choice
of basis as above. Quotienting by GIT with the natural action of 
$Sl(V)$ on $R$, induced from its natural action on $\mathbb{P}$, 
we obtain the wanted
projective coarse moduli space.

\begin{definition}
\label{4}
Let $G=O(r,\mathbb{C})$ or $Sp(r,\mathbb{C})$. A principal $G$-sheaf on $X$ is a
tensor $\varphi :E\otimes E\longrightarrow \SO _{X}$ symmetric or
antisymmetric which induces an isomorphism $E_{|U}\longrightarrow
E_{|U}^{\ast }$ on the open set $U$ where $E$ is locally free. We call it
(semi)stable if for all isotropic subsheaves $F\subseteq E$ it is 
\[
\chi _{F}+\chi _{F^{\perp }}(\leq )\chi _{E} 
\]
\end{definition}

\begin{theorem}
\label{5}
For any positive polynomial $\delta $ of degree at most $n-1$, 
a principal $G$-sheaf on $X$ $(G=O(r,\mathbb{C))}$ or
$Sp(r,\mathbb{C})$ is $\delta $-(semi)stable if and only if 
it is (semi) stable, so there is a coarse
projective moduli space of stable-equivalence classes of semistable 
principal $G$-sheaves.
\end{theorem}

\textbf{The remaining classical group }$G=SO(r,\CC)$. 
Define a principal $SO(r,\CC)$-sheaf to be a triple
$(E,\varphi,\psi)$, where $(E,\varphi)$ is a principal 
$O(r,\CC)$-sheaf and $\varphi$ is an isomorphism between
$\det(E)$ and $\SO_X$ such that $\det(\varphi)=\psi^2$.
Note that for each $O(r,\CC)$-sheaf $(E,\varphi)$, there is
at most two distinct $SO(r,\CC)$-sheaves, namely
$(E,\varphi,\psi)$ and $(E,\varphi,-\psi)$. If $\rk(E)$ is odd,
these two objects are isomorphic. This is why for $SO(2m+1,\CC)$
we can forget the third datum $\psi$. But if $\rk(E)$ is even,
these two objects might not be isomorphic.
With the same definition of
(semi)stability as in Definition \ref{4}, Theorem \ref{5} 
still holds in this case (i.e. the
added datum does not alter the GIT notion of stability) 
so \textit{we obtain a coarse
projective moduli space in the case $G$ is any classical group.}

\section{Principal sheaves on a reductive group}
\label{sec4}

Tensors considered in Section \ref{sec3} 
were all $(s,0)$ tensors, but with the same
machinery we could have worked with (semi)stability and coarse projective
moduli space of $(s,1)$-tensors. In particular we need in 
this section $(2,1)$-tensors $\varphi :E\otimes E\longrightarrow
E^{\ast \ast}$, 
for which $\delta $-(semi)stability is defined by the fact that 
for all weighted filtration $(E_{1}\subset ...\subset E_{t},
m_{1},...,m_{t}>0)$ of $E$, it is 
\[
\sum m_{i}(r\chi _{E_{i}}-r_{i}\chi _{E})+\delta\; \mu (E.,m.,\varphi )
\; (\leq )\;0
\]
where 
\[
\mu =\min \left\{ \lambda _{i}+\lambda _{j}-\lambda _{k}|0\neq 
\overline{\varphi }:E_{i}\otimes E_{j}\longrightarrow
E^{\ast\ast}/E_{k-1}^{\ast\ast}
\right\} 
\]
For fixed value of rank and Chern classes, there is a projective coarse
moduli space of stable equivalence classes 
of $\delta $-semistable $(2,1)$ tensors
on $X.$

\begin{definition}
Let $X$ be a projective variety, and $G$ an algebraic group. A 
principal $G$-sheaf $\SP$ is a triple $(E,\varphi ,\xi )$ where
$(E,\varphi)$ is a $(2,1)$-tensor on $X$
$$
\varphi :E\otimes E\longrightarrow E^{\ast\ast}
$$ 
such that for the 
points $x$ of the open set $U_{E}$ where $E$ is locally free, 
$\varphi (x)$ is isomorphic to the
structure tensor $\varphi _{\fg'}:\fg'\otimes 
\fg'\longrightarrow \fg'$ of the Lie algebra $\fg'$ tangent to 
the commutator subgroup $G^{\prime }=[G,G]$ (in particular,
there is an associated $Aut(\fg')$-bundle $P_{U_E}$ on $U_E$),
and $\xi$ is a reduction of $P_{U_E}$
to $G$, via  $Ad:G\longrightarrow Aut(\fg')$.
\end{definition}

Obviously, if $E$ is locally free, we recover the usual notion of 
principal $G$-bundle.

\begin{definition}
\label{7}
Let $G$ be a connected reductive group. We say that a principal 
$G$-sheaf $\mathcal{P}=(E,\varphi ,\xi )$ is semistable if $E$ 
is semistable. We say it
is strictly semistable if there is not a Lie algebra filtration, i.e. 
\[
E_{1}\subseteq ...\subseteq E_{t}\subseteq E_{t+1}=E\; 
\]
such that $[E_{i},E_{j}]\subseteq E_{i+j}^{\ast\ast}$ with all 
\[
\frac{\chi _{E_{i}}}{r_{i}}=\frac{\chi _{E}}{r} 
\]
\end{definition}

\begin{theorem}
There is a projective coarse moduli stable of equivalence classes of
semistable principal $G$-sheaves on $X$ of fixed topological type.
\end{theorem}

\textbf{Comment on the proof}: It is a very long proof, 
parallel to the proof of Ramanathan \cite{R2},
which will appear
published elsewhere. 
Because of the
nondegeneracy of the Killing form of the semisimple Lie algebra 
$\fg'$, the factor $\mu (E.,m.,\varphi )$ is always nonpositive, so for a
polynomial $\delta $ of degree zero, and small with respect to
the invariants 
of $X$, $E$, our notion of (semi)stability of $\mathcal{P}=
(E,\varphi ,\xi )$ is
equivalent to the $\delta $-(semi)stability of the $(2,1)$-tensor 
$(E,\varphi)$. 
It does not assure the existence of a moduli space, 
because it must also be checked
that the extra datum of reduction $\xi $ does not alter the (semi)stability
in the sense of GIT of the corresponding point of the 
$Sl(V)$-acted 
projective space, which is the main bulk of the proof.

Finally, we need some considerations on root spaces in order to re-state
(semi)stability in a way which is equivalent, but 
more convenient to check that it coincides with 
Ramanathan's (semi) stability when $\dim X=1$. 
Recall
from \cite{Bo} 
that a $\ft$-root decomposition 
\[
\fg'=\bigoplus_ {\alpha \in R_{t}\cup \{0\}}\fg^{\prime\alpha }
\]
of the Lie algebra $\fg'$ arises whenever an abelian algebra 
$\ft\subseteq \fg'$ is given, 
not necessarily a Cartan algebra, in particular for
the center $\ft=\fz(\fl(\fh'))$ of the Levi component 
$\fl(\fh')$ of any parabolic
subalgebra $\fh'\subset \fg'$. 
In this case a system of simple $\ft$-roots 
(or decomposition $R_{\ft}=R_{\ft}^{+}\cup R_{\ft}^-$) 
is naturally given, so the set $R_{\ft}\cup \{0\}$
has a natural partial ordering ($\alpha \leq \beta $ if $\beta $ is the sum
of $\alpha $ with a sum of simple $\ft$-roots). Denote 
$\fg'_{(\leq \alpha )}=
\oplus_{\beta \leq \alpha } \fg^{\prime\beta }$ and 
analogously $\fg'_{(<\alpha )}$. We also write $R_{\fh'}$ for 
$R_{\ft}$. Both are invariant by
the adjoint action of $\fh'$, thus by the inner automorphism action of the
corresponding parabolic subgroup $H'$ of the group 
$G'$, so the
analogous subalgebras $\fg_{(\leq \alpha )}^{\prime }$ and 
$\fg_{(<\alpha)}^{\prime }$ of the Lie algebra $\fg'$ are also 
$H'$-invariant.

Let $\mathcal{P}=(E,\varphi ,\xi )$ be a principal sheaf on $X$,
 having 
on $U_{E}$ a further $H\hookrightarrow G$ reduction to a parabolic 
subgroup $H$, let $H'=H\cap G'$, 
and let $\alpha \in R_{\fh'}\cup \{0\}$ where
 $\fh'=\operatorname{Lie}(H')$ as before. We define 
$E_{(\leq \alpha )}$ and $E_{(<\alpha )}$ as the saturated extensions 
to $X$ of the vector bundles on $U_{E}$ associated to this reduction 
by the above representation of $H'$ on $\fg_{(\leq \alpha )}^{\prime }$ 
and $\fg_{(<\alpha )}^{\prime }$, and define 
$E^{\alpha }$ as $E_{(\leq \alpha )}/E_{(<\alpha )}.$

\begin{proposition}
\label{9}
A semistable principal $G$-sheaf $\mathcal{P}=(E,\varphi ,\xi )$ on $X$
is stable if and only if for all reductions of 
$\mathcal{P}|_{U_{E}}$ to a maximal
parabolic subgroup $H$ and $\alpha \in R_{\fh'}\cup \{0\}$, it is 
\[
\frac{\chi_{E^{\alpha }}}{\rk(E^{\alpha })}<\frac{\chi_E}
{\rk(E)} 
\]
\end{proposition}

\medskip

\textbf{The case dim X=1.} 
In this case we have $U_E=X$, then a principal sheaf is equivalent
to a principal bundle. The polynomials 
$\chi_{E^{\alpha}}$ 
and $\chi_E$ in Proposition \ref{9} can be replaced by 
$\deg E^{\alpha }$, $\deg E$, and,
being $\deg E=0$, the strict semistability of a semistable principal bundle
amounts to $\deg E^{\alpha }<0$ for all $\alpha $ as in Proposition 
\ref{9}. A short argument shows the following

\begin{proposition}
If $\dim X=1$, a principal bundle  
$\mathcal{P}=(E,\varphi ,\xi )$ is (semi)stable (Definition \ref{7})
if and only if for all reductions $P(H)$ to a maximal
parabolic subgroup $H$ of $G$, we have 
$\deg(P(H,\fh))(\leq)0$,
where $P(H,\fh)$ is the vector bundle associated to $P(H)$
by the adjoint representation of $H$ in its Lie algebra $\fh$.
\end{proposition}

\bigskip Therefore, in  case $\dim X =1$, we obtain exactly that our
notion of (semi)stability coincides with Ramanathan's definition 
\cite[Remark 2.2]{R1}.

\begin{remark}
In the case of $G$ classical group, it would be interesting
to compare the definitions of sections \ref{sec3} and \ref{sec4}.
If $\dim X=1$, then they coincide.
\end{remark}

\noindent Tom\'as L. G\'omez: School of Mathematics, Tata Institute of 
Fundamental Research, Bombay 400 005 (India)

\noindent \texttt{tomas@math.tifr.res.in}

\medskip

\noindent Ignacio Sols: Departamento de Algebra, 
Facultad de Matem\'aticas,
Universidad Complutense de Madrid, 28040 Madrid (Spain)

\noindent \texttt{sols@eucmax.sim.ucm.es}

\end{document}